\documentclass{article}
\usepackage{amsmath,amssymb,verbatim}

\usepackage{times}

\usepackage[dvips]{graphicx}
\usepackage{amsfonts}
\usepackage{amsthm}
\usepackage{epsfig}

\theoremstyle{plain}
\newtheorem{theorem}{Theorem}[section]
\newtheorem{lemma}[theorem]{Lemma}
\newtheorem{proposition}[theorem]{Proposition}
\newtheorem{corollary}[theorem]{Corollary}
\theoremstyle{definition}
\newtheorem{definition}[theorem]{Definition}

\newtheorem{remark}[theorem]{Remark}

\newtheorem{example}[theorem]{Example}

\newcommand{\xs}{x_1,\ldots,x_n}                
\newcommand{\Fs}{F_1,\ldots,F_q}                

\newcommand{\HH}{{\mathcal{H}}}                  
\newcommand{\D}{\Delta}                         
\newcommand{\al}{\alpha}                         
\newcommand{\be}{\beta}                          
              
\newcommand{\rmf}[1]{\setminus \{#1\}}          
\newcommand{\void}[1]{}

\newcommand{\cocoa}{\mbox{\rm C\kern-.13em o\kern-.07 em C\kern-.13em o\kern-.15em A}} 
\newcommand{\cocoax}{\mbox{C\kern-.13em o\kern-.07 em C\kern-.13em o\kern-.15em A}} 
\newcommand{\cocoal}{\mbox{\rm C\kern-.13em o\kern-.07 em C\kern-.13em o\kern-.15emA\kern-.1em L}}

\newcommand{\todo}[1]{\vspace{5 mm}\par \noindent
\marginpar{\textsc{ToDo}}
\framebox{\begin{minipage}[c]{0.95 \textwidth}
\tt #1 \end{minipage}}\vspace{5 mm}\par}

\renewcommand{\todo}[1]{}

\newcommand{\idiot}[1]{\vspace{5 mm}\par \noindent
\framebox{\begin{minipage}[c]{0.95 \textwidth}
\tt #1 \end{minipage}}\vspace{5 mm}\par}

\renewcommand{\idiot}[1]{}

\newcommand{\seq}{\subseteq}    
\newcommand{\es}{\emptyset}             
\newcommand{\s}[1]{\{#1\}}              

\renewcommand{\leq}{\leqslant}          
\renewcommand{\geq}{\geqslant}  

\date{\today}

\author{Massimo Caboara\thanks{Department of Mathematics, University
of Pisa, caboara@dm.unipi.it.} \and 
Sara Faridi\thanks{Department of
Mathematics and Statistics, Dalhousie University, Halifax, Canada, 
faridi@mathstat.dal.ca. 
Research supported by NSERC.}}

\title{Odd-Cycle-Free Facet Complexes and the K\"onig property}
\begin{document}

\maketitle
\begin{abstract} We use the definition of a simplicial cycle to define
 an odd-cycle-free facet complex (hypergraph). These are facet
 complexes that do not contain any cycles of odd length. We show that
 besides one class of such facet complexes, all of them satisfy the
 K\"onig property. This new family of complexes includes the family of
 balanced hypergraphs, which are known to satisfy the K\"onig
 property. These facet complexes are, however, not Mengerian; we give
 an example to demonstrate this fact.
 \end{abstract}

\section{Introduction} 
 
 Simplicial trees were introduced by the second author in \cite{F1} in
 order to generalize algebraic structures based on graph trees. More
 specifically, the facet ideal of a simplicial tree, which is the
 ideal generated by the products of the vertices of each facet of the
 complex in the polynomial ring whose variables are the vertices of
 the complex, is a normal ideal (\cite{F1}), is always sequentially
 Cohen-Macaulay (\cite{F2}) and one can determine exactly when the
 quotient of this ideal is Cohen-Macaulay based on the combinatorial
 structure of the tree (\cite{F3}). These algebraic results that
 generalize those associated to simple graphs, and are intimately tied
 to the combinatorics of the simplicial complex, have suggested that
 this is a promising definition of a tree in higher dimension.  This
 fact was most recently confirmed when the authors, while searching
 for an efficient algorithm to determine when a given complex is a
 tree, produced a precise combinatorial description for a simplicial
 cycle that has striking resemblance to that of a graph cycle
 (\cite{CFS}). The main idea here is that a complex (or a simple
 hypergraph) is a tree if and only if it does not contain any
 ``holes'', or any cones over holes.  Moreover, our definition of a
 simplicial cycle, though more restrictive than that defined for
 hypergraphs by Berge~\cite{B1,B2}, satisfies the hypergraph definition as
 well. In a way, simplicial cycles are ``minimal'' hypergraph cycles,
 in the sense that once a facet is removed, what remains is not a
 cycle anymore, and does not contain one.

 Once the concept of a ``minimal'' cycle is in place, a natural
 question that arises is whether the length of such a cycle bears any
 meaning in terms of properties of the complex? In graph theory
 bipartite graphs are characterized as those that do not contain any
 odd cycles. One of their strongest features is that they satisfy the
 K\"onig property. Our purpose in this paper is to investigate whether
 simplicial complexes (or hypergraphs) not containing odd simplicial
 cycles, which we call \emph{odd-cycle-free} complexes, also satisfy
 this property. It turns out that besides one family, all
 odd-cycle-free complexes do satisfy the K\"onig property
 (Theorem~\ref{t:mainthm}).  The proof uses tools from hypergraph
 theory, as well as Berge's recently proved Strong Perfect Graph
 Conjecture (\cite{C,CRST}).

 A much more general notion of a cycle already exists in hypergraph
 theory (\cite{B1,B2}); we call these \emph{hyper-cycles}
 (Definition~\ref{d:hyper-cycle}) to avoid confusion. It is known that
 hypergraphs that do not contain odd hyper-cycles are \emph{balanced},
 and hence satisfy the K\"onig property.  The class of odd-cycle-free
 complexes which we study in this paper includes the class of simple
 hypergraphs that do not contain odd hyper-cycles, and hence our
 results generalize those already known for hypergraphs.  We discuss
 these inclusions in Section~\ref{s:hyper-cycles}.

 Simis, Vasconcelos and Villarreal showed in~\cite{SVV} that facet
 ideals of bipartite graphs are normally torsion free, and hence
 normal. Recently Herzog, Hibi, Trung and Zheng~\cite{HHTZ} have
 generalized their result and shown that facet ideals of Mengerian
 complexes (hypergraphs) are normally torsion free. This includes the
 class of simple hypergraphs that do not contain odd hyper-cycles, and
 more generally, balanced hypergraphs. We demonstrate in
 Section~\ref{s:hyper-cycles} that odd-cycle-free complexes are not
 necessarily Mengerian, and hence their facet ideals are not
 necessarily normally torsion-free, although they could still be
 normal ideals.
 
 While this paper refers to simplicial or facet complexes most of the
 time for the statements, it is important to know that these
 structures, for our purposes, are essentially simple hypergraphs. The
 original work on higher dimensional trees and cycles was done in the
 context of commutative algebra, where a rich tradition of studying
 ideals associated to simplicial complexes was already in place.  This
 paper, on the other hand, uses many results from hypergraph
 theory. For this reason, and for the sake of consistency, in the
 introductory parts of the paper, we give a careful review of all the
 structures that we use and demonstrate how it is possible to move
 between complexes and hypergraphs without losing the validity of any
 of our statements.

\section{Facet complexes, trees, and cycles}\label{s:first-section} 

We define the basic notions related to facet complexes. More details and
examples can be found in~\cite{F1,F3}.

\begin{definition}[Simplicial complex, facet]
A \emph{simplicial complex} $\D$ over a finite set of vertices $V$
is a collection of subsets of $V$, with the property that if $F \in
\D$ then all subsets of $F$ are also in $\D$. An element of $\D$ is
called a \emph{face} of $\D$, and the maximal faces are called
\emph{facets} of $\D$.
\end{definition}

Since we are usually only interested in the facets, rather than all
faces, of a simplicial complex, it will be convenient to work with the
following definition:

\begin{definition}[Facet complex]
  A \emph{facet complex} over a finite set of vertices $V$ is a set
  $\D$ of subsets of $V$, such that for all $F,G\in\D$, $F\seq G$
  implies $F=G$.  Each $F\in\D$ is called a \emph{facet} of $\D$.
\end{definition}

\begin{remark}[Equivalence of simplicial complexes and facet complexes]
  The set of facets of a simplicial complex forms a facet
  complex. Conversely, the set of subsets of the facets of a facet
  complex is a simplicial complex. This defines a one-to-one
  correspondence between simplicial complexes and facet complexes. 
  In this paper, we will work primarily with facet complexes.
\end{remark}

 We now generalize some notions from graph theory to facet complexes. Note
 that a graph can be regarded as a special kind of facet complex,
 namely one in which each facet has cardinality 2.

\begin{definition}[Path, connected facet complex]\label{def:connected}
  Let $\D$ be 
  a facet complex.  A sequence of facets $F_1,\ldots,F_n$ is
  called a \emph{path} if for all $i=1,\ldots,n-1$, $F_i\cap
  F_{i+1}\neq\es$.  We say that two facets $F$ and $G$ are {\em
  connected} in $\D$ if there exists a path $F_1,\ldots,F_n$ with
  $F_1=F$ and $F_n=G$. Finally, we say that $\D$ is
  \emph{connected} if every pair of facets is connected.
\end{definition}

In order to define a tree, we borrow the concept of \emph{leaf} from
graph theory, with a small change.

\begin{definition}[Leaf, joint]\label{d:leaf} Let $F$ be a
 facet of a facet complex $\D$.  Then $F$ is called a \emph{leaf} of
$\D$ if either $F$ is the only facet of $\D$, or else there exists
some $G \in \D \rmf{F}$ such that for all $H\in\D\rmf{F}$, we have $H\cap
F \subseteq G$. The facet $G$ above is called a \emph{joint} of the leaf
$F$ if $F \cap G \neq \emptyset$.
\end{definition}

It follows immediately from the definition that every leaf $F$
contains at least one {\em free vertex}, i.e., a vertex that belongs
to no other facet.

\begin{example}\label{example11} 
In the facet complex $\D=\s{F, G, H}$, $F$ and $H$ are leaves, but
 $G$ is not a leaf.  Similarly, in $\D'=\s{A,B,C}$, the only leaves
 are $A$ and $C$.
\[
\D=\begin{tabular}{c}\epsfig{file=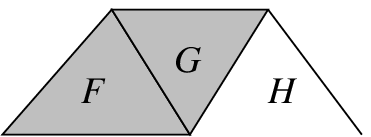, height=.6in}\end{tabular}
\qquad
\D'=\begin{tabular}{c}\epsfig{file=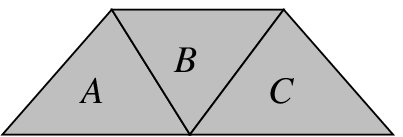, height=.6in}\end{tabular}
\]
\end{example}

In Example~\ref{example11} as well as in the rest of this paper, we
use a shaded $n$-polygon to display a facet with $n$ vertices. So we
can think of the facet complex $\D$ in the example above as if the
vertices were labeled with $x,y,z,u$, such that $F=\{x,y,z\}$,
$G=\{y,z,u\}$ and $H=\{u,v\}$.

\begin{definition}[Forest, tree]\label{d:tree} 
  A facet complex $\D$ is a \emph{forest} if every nonempty subset of
  $\D$ has a leaf. A connected forest is called a \emph{tree} (or
  sometimes a {\em simplicial tree} to distinguish it from a tree in
  the graph-theoretic sense).
\end{definition}

It is clear that any facet complex of cardinality one or two is a forest.
When $\D$ is a graph, the notion of a simplicial tree coincides with
that of a graph-theoretic tree.

\begin{example}\label{e:free-example} 
  The facet complexes in Example~\ref{example11} are trees. The facet
  complex pictured below has three leaves $F_1$, $F_2$ and $F_3$;
  however, it is not a tree, because if one removes the facet $F_4$,
  the remaining facet complex has no leaves.
\begin{center}
 \epsfig{file=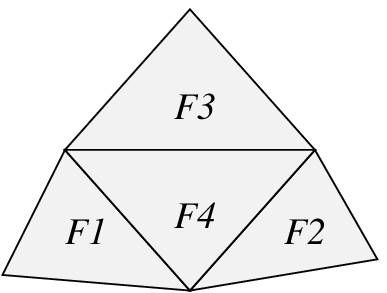, height=1in}
\end{center}
\end{example}

\begin{definition}[Minimal vertex cover, Vertex covering number] 
  Let $\Delta$ be a facet complex with vertex set $V$ and facets
  $\Fs$. A \emph{vertex cover} for $\Delta$ is a subset $A$ of $V$,
  with the property that for every facet $F_i$ there is a vertex $v
  \in A$ such that $v \in F_i$. A \emph{minimal vertex cover} of
  $\Delta$ is a subset $A$ of $V$ such that $A$ is a vertex cover, and
  no proper subset of $A$ is a vertex cover for $\Delta$. The smallest
  cardinality of a vertex cover of $\Delta$ is called the \emph{vertex
    covering number} of $\Delta$ and is denoted by $\al(\D)$.
\end{definition}

\begin{definition}[Independent set, Independence number]
  Let $\D$ be a facet complex. A set $\{F_1, \ldots, F_u\}$ of
  facets of $\D$ is called an \emph{independent set} if $F_i \cap F_j
  =\emptyset$ whenever $i \neq j$. The maximum possible cardinality of
  an independent set of facets in $\D$, denoted by $\beta(\D)$, is
  called the \emph{independence number} of $\D$. An independent set of
  facets which is not a proper subset of any other independent set is
  called a \emph{maximal independent set} of facets.
\end{definition}

\subsection{Cycles}\label{s:cycles}

In this section, we define a simplicial cycle as a minimal facet
complex without leaf. This in turn characterizes a tree as a connected
cycle-free facet complex. The main point is that higher-dimensional
cycles, like graph cycles, possess a particularly simple structure:
each cycle is either equivalent to a ``circle'' of facets with
disjoint intersections, or to a cone over such a circle.

\begin{definition}[Cycle]\label{d:cycle}
A nonempty facet complex $\D$ is called a {\em cycle} (or a
\emph{simplicial cycle}) if $\D$ has no leaf but every nonempty proper
subset of $\D$ has a leaf.
\end{definition}

Equivalently, $\D$ is a cycle if $\D$ is not a forest, but every
proper subset of $\D$ is a forest.  If $\D$ is a graph,
Definition~\ref{d:cycle} coincides with the graph-theoretic definition
of a cycle. The next remark is an  immediate consequence of the
definitions of cycle and forest.

\begin{remark}[A forest is a cycle-free facet complex]
  A facet complex is a forest if and only if it does not contain a
  cycle.
\end{remark}

We now provide a complete characterization of the structure of cycles
as described in~\cite{CFS}.

\begin{definition}[Strong neighbor]\label{d:strong-neighbor} Let $\D$ 
be a facet complex and $F,G\in\D$. We say that $F$ and $G$ are
\emph{strong neighbors}, written $F\sim_{\D} G$, if $F\neq G$ and for
all $H\in\D$, $F\cap G\seq H$ implies $H=F$ or $H=G$.
\end{definition}

  The relation $\sim_{\D}$ is symmetric, i.e., $F\sim_{\D} G$ if and only if
  $G\sim_{\D} F$. Note that if $\D$ has
  more than two facets, then $F\sim_{\D} G$ implies that $F \cap G \neq
  \emptyset$.

\begin{example} For the facet complex $\D'$ in Example~\ref{example11}, 
$A \not \sim_{\D'}C$, as their intersection lies in the facet
$B$. However, $B\sim_{\D'}C$ and similarly $B\sim_{\D'}A$.
\end{example}

A cycle can be described as a sequence of strong neighbors. 

\begin{theorem}[Structure of a cycle~(\cite{CFS})]\label{t:cycle-structure}
 Let $\D$ be a facet complex.  Then $\D$ is a cycle if and only if the
facets of $\D$ can be written as a sequence of strong neighbors
$F_1\sim_{\D} F_2 \sim_{\D} \ldots \sim_{\D} F_n \sim_{\D} F_1$ such
that $n\geq 3$, and for all $i,j$ $$F_i \cap F_j = \bigcap_{k=1}^n F_k
\quad\mbox{ if } j \neq i-1, i,i+1\ (\mbox{\rm mod } n).$$
\end{theorem}

The implication of Theorem~\ref{t:cycle-structure} is that a
 simplicial cycle has a very intuitive structure: it is either a
 sequence of facets joined together to form a circle (or a ``hole'')
 in such a way that all intersections are pairwise disjoint (this is
 the case where the intersection of all the facets is the empty set in
 Theorem~\ref{t:cycle-structure}), or it is a cone over such a
 structure.

\begin{example}\label{e:cycle-structure} The facet complex $\D$ is a cycle.
  The facet complex $\Gamma$ is a cycle and is also a cone over the
  cycle $\Gamma'$.

\begin{center}
$\D=\begin{tabular}{c}\epsfig{file=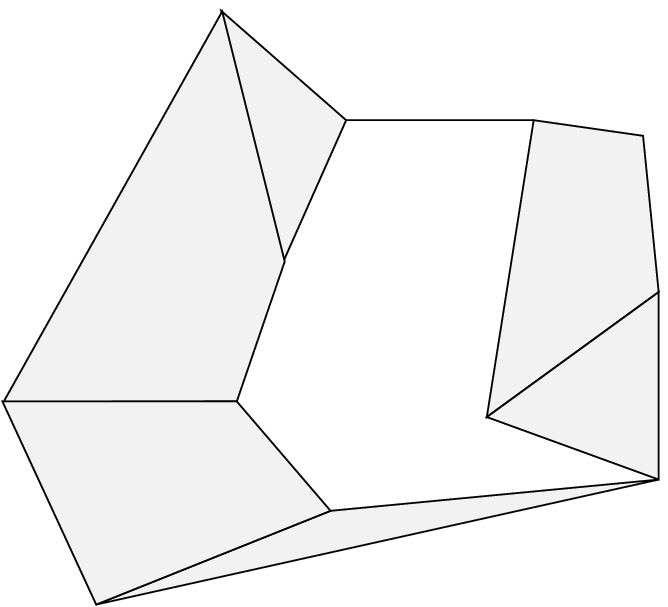, width=1.2in}\end{tabular}
\qquad
\Gamma=\begin{tabular}{c}\epsfig{file=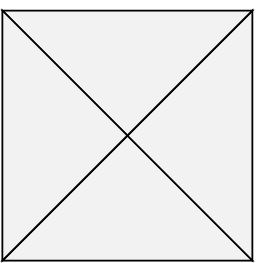, height=.7in}\end{tabular}
\qquad
\Gamma'=\begin{tabular}{c}\epsfig{file=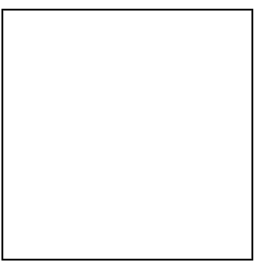, height=.7in}\end{tabular}$
\end{center}
\end{example}

The next example demonstrates the impact of the second condition of
being a cycle in Theorem~\ref{t:cycle-structure}.

\begin{example}\label{e:cycle-is} The facet complex $\D$ has no leaves 
but is not a cycle, as its proper subset $\D'$ (which is indeed a
 cycle) has no leaves. However, we have
 $F_1\sim_{\D}F_2\sim_{\D}G\sim_{\D}F_3\sim_{\D}F_4\sim_{\D}F_1$, and
 these are the only pairings of strong neighbors in $\D$.

\[
\D=\begin{tabular}{c}\epsfig{file=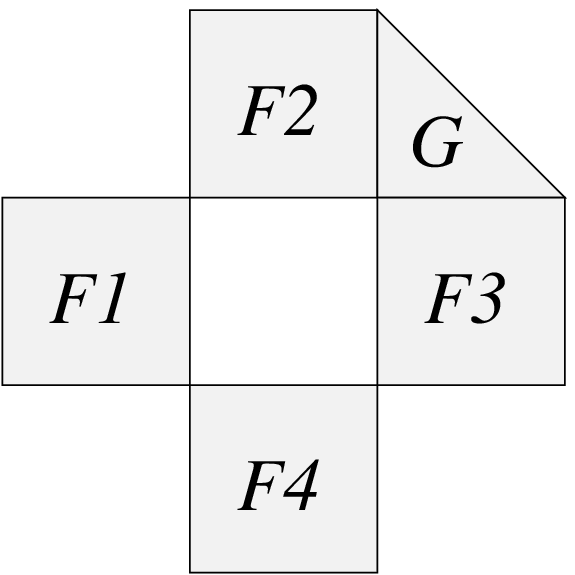, height=1in}\end{tabular}
\qquad
\D'=\begin{tabular}{c}\epsfig{file=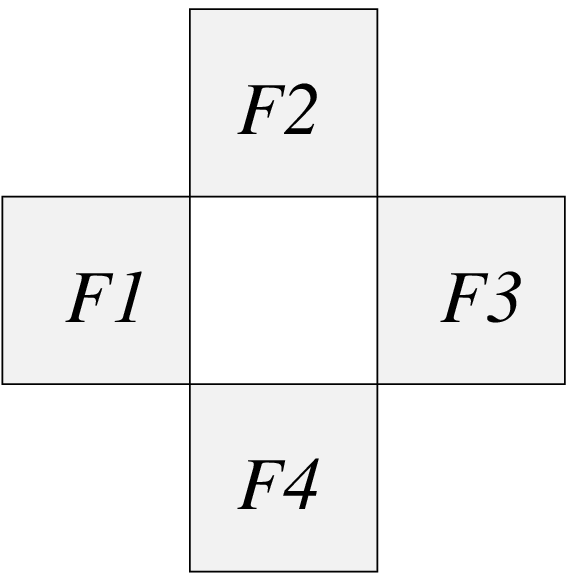, height=1in}\end{tabular}
\]

\end{example}

A property of cycles that we shall use often in this paper is the following.

\begin{lemma}\label{l:3-cycle} Let $F_1,F_2,F_3$ be facets of a 
facet complex $\D$, such that $F_i\cap F_j\neq \emptyset$ for $i,j \in
\{1,2,3\}$, and $F_1\cap F_2 \cap F_3 =\emptyset$. Then
$\Gamma=\{F_1,F_2,F_3\}$ is a cycle.
\end{lemma}

     \begin{proof} Since $\Gamma$ has three facets, all its proper subsets 
     are forests. So if $\Gamma$ is not a cycle, then it must contain
     a leaf.  Say $F_1$ is a leaf, and $F_2$ is its joint. So we have
     $\emptyset \neq F_1 \cap F_3 \subseteq F_2$, which implies that
     $F_1\cap F_2 \cap F_3 \neq \emptyset$; a contradiction.
     \end{proof}

\section{Facet complexes as simple hypergraphs}

\subsection{Graph theory terminology} 

\begin{definition}[Induced subgraph]  Let $G$ be a graph with vertex set 
$V$. A subgraph $H$ of $G$ with vertex set $W\subseteq V$ is called an
\emph{induced subgraph} of $G$ if for each $x,y \in W$, $x$ and $y$
are connected by an edge in $H$ if and only if they are connected by
an edge in $G$.
\end{definition}

\begin{definition}[Clique of a graph] A clique of a graph $G$ is a   
complete subgraph of $G$; in other words a subgraph of $G$ whose every
two vertices are connected by an edge. \end{definition}

\begin{definition} [Chromatic number] The \emph{chromatic number} of 
a graph $G$ is the smallest number of colors needed to color the
 vertices of $G$ so that no two adjacent vertices (vertices that
 belong to the same edge) share the same color.
\end{definition}

\begin{definition}[Complement of a graph] The complement of a graph $G$, 
denoted by $\overline{G}$, is a graph over the same vertex set as $G$
whose edges connect non-adjacent vertices of $G$.
\end{definition}

\begin{definition}[Perfect graph] A graph $G$ is \emph{perfect}
 if for every induced subgraph $G'$ of $G$, the chromatic number of
$G'$ is equal to the size of the largest clique of $G'$.
\end{definition}

We call $G$ a \emph{minimal imperfect} graph if it is not perfect but
all proper induced subgraphs of $G$ are perfect. There is a
characterization of minimal imperfect graphs that was conjectured by
Berge and known for a long time as ``Strong Perfect Graph
Conjecture'', and was proved recently by Chudnovsky, Robertson,
Seymour and Thomas~\cite{CRST}; see also~\cite{C}.

\begin{theorem}[Strong Perfect Graph Theorem (\cite{CRST})]\label{t:SPGT} 
The only minimal imperfect graphs are odd cycles of length $\geq 5$
 and their complements.
\end{theorem}

\subsection{Hypergraphs}

A hypergraph is simply a higher dimensional graph.

\begin{definition}[Hypergraph, simple hypergraph (\cite{B1})] 
Let $V=\{\xs\}$ be a finite set.  A \emph{hypergraph} on $V$ is a
family $\HH=(F_1,\ldots,F_m)$ of subsets of $V$ such that
\begin{enumerate}
\item $F_i\neq \emptyset$ for $i=1,\ldots,m$;
\item $\displaystyle V=\cup_{i=1}^m F_i$.
\end{enumerate}
Each $F_i$ is called an \emph{edge} of $\HH'$.  If, additionally, we
have the condition: $F_i \subset F_j \Longrightarrow i=j$, then $\HH$
is called a \emph{simple hypergraph}.
\end{definition}

A graph is a hypergraph in which an edge consists of exactly two vertices.

\begin{definition}[Partial hypergraph] A \emph{partial hypergraph} of a 
hypergraph $\HH=\{F_1,\ldots,F_m\}$ is a subset $\HH'=\{F_j\ |\ j \in
J\}$, where $J\subseteq \{1,\ldots,m\}$.
\end{definition}

It is clear that a facet complex $\D$ is a simple hypergraph on its
set of vertices, and a partial hypergraph is just a subset of
$\D$. For this reason, we are able to borrow the following definitions
from hypergraph theory. The main source for these concepts is Berge's
book~\cite{B1}.

\begin{definition}[Line graph of a hypergraph] Given a hypergraph 
$\HH=\{F_1,\ldots,F_m\}$ on vertex set $V$, its \emph{line graph}
$L(\HH)$ is a graph whose vertices are points $e_1,\ldots,e_m$
representing the edges of $\HH$, and two vertices $e_i$ and $e_j$ are
connected by an edge if and only if $F_i \cap F_j \neq \emptyset$.
\end{definition}

\begin{definition}[Normal hypergraph (\cite{L})] A hypergraph 
$\HH$ with vertex set $V$ is \emph{normal} if every partial hypergraph
$\HH'$ satisfies the \emph{colored edge property}, i.e.
$q(\HH')=\delta(\HH')$, where

\begin{itemize}

\item $q(\HH') =$ \emph{chromatic index of $\HH'$}, which is the
 minimum number of colors required color the edges of $\HH'$ in such a
 way that two intersecting edges have different colors; and

\item $\delta(\HH')= {\rm max}_{x\in V} \{\mbox{number of edges of $\HH'$ that 
contain } x\}$.
\end{itemize}
Clearly, we always have $q(\HH')\geq \delta(\HH')$.
\end{definition}

\begin{definition}[Helly property]
Let $\HH=\{\Fs\}$ be a simple hypergraph, or equivalently, a facet
complex. Then $\HH$ is said to satisfy the \emph{Helly property} if
every intersecting family of $\HH$ is a star; i.e., for every $J
\subseteq \{1,\ldots,q\}$ $$F_i\cap F_j \neq \emptyset \mbox{ for all
} i,j \in J \Longrightarrow \bigcap_{j\in J} F_j \neq \emptyset.$$
\end{definition}

From the above definitions, the following statement, which we shall
rely on for the rest of this paper, makes sense.

\begin{theorem}[\cite{B1} page 197]\label{t:Helly-perfect-normal} 
A simple hypergraph (or facet complex) $\HH$ is normal if and only if
$\HH$ satisfies the Helly property and $L(\HH)$ is a perfect graph.
\end{theorem}
 
\section{Odd-Cycle-Free complexes}

As we discussed in the previous section, a facet complex is a simple
hypergraph. One particular property of facet complexes that we are
interested in is the K\"onig property.

\begin{definition}[K\"onig property] A facet complex $\D$ satisfies the 
\emph{K\"onig property} if $\al(\D)=\be(\D)$.\end{definition}

\begin{definition}[Odd-Cycle-Free complex] We call a facet complex 
\emph{odd-cycle-free} if it contains no cycles of odd
length.\end{definition}

It is well known that odd-cycle-free graphs, which are known as
bipartite graphs, satisfy the K\"onig property. In higher dimensions,
this property is enjoyed by simplicial trees~\cite{F3}, and complexes
that do not contain special odd cycles, also known as \emph{balanced
hypergraphs}~\cite{B1,B2}.  The class of odd-cycle-free complexes includes
all such complexes (see Section~\ref{s:hyper-cycles}).

It is therefore natural to ask if odd-cycle-free complexes satisfy the
K\"onig property. The answer to this question is mostly positive:
besides one specific class of odd-cycle-free complexes, all of them do
satisfy the K\"onig property.

We begin with from the following fact due to Lov\'asz~\cite{L} (see
also~\cite{B1} page 195).

\begin{theorem}[Normal hypergraphs satisfy K\"onig]\label{t:normal-konig} 
The hypergraph $\HH$ is normal if and only if every partial hypergraph
of $\HH$ satisfies the K\"onig property.
\end{theorem}

 We can hence prove that a facet complex $\D$ (and its subsets)
satisfy the K\"onig property by showing that $\D$ is normal.

\begin{theorem}[Odd-Cycle-Free complexes that are normal]\label{t:odd-cycle-free-normal}  
If $\D$ is a facet complex that is odd-cycle-free and $L(\D)$ does not
contain the complement of a 7-cycle as an induced subgraph, then $\D$
is normal.
\end{theorem}

By Theorem~\ref{t:Helly-perfect-normal}, it suffices to show
that $\D$ satisfies the Helly property and $L(\D)$ is perfect.  We
show these two properties separately.

\begin{proposition}[Odd-Cycle-Free complexes satisfy Helly property]\label{p:odd-cycle-free-Helly} If $\D$ is an odd-cycle-free
 facet complex, then it satisfies the Helly property.
\end{proposition}

     \begin{proof} Suppose $\D$ does not satisfy the Helly property, 
      so it contains an intersecting family that is not a star. In
      other works, there exists $\Gamma=\{F_1,\ldots,F_m\} \subseteq
      \D$ such that $$F_i \cap F_j \neq \emptyset \mbox{ for }  i,j \in
      \{1,\ldots,m\}, \mbox{ but } \bigcap_{j=1}^m F_j =\emptyset.$$ 

       We use induction on $m$. If $m=3$, from Lemma~\ref{l:3-cycle}
       it follows that $\Gamma$ is a 3-cycle. 

       Suppose now that $m>3$ and we know that every intersecting
       family of less than $m$ facets that is not a star contains an
       odd cycle. Let $\Gamma$ be an intersecting family of $m$ facets
       $F_1,\ldots,F_m$, such that every $m-1$ facets of $\Gamma$
       intersect (otherwise by the induction hypothesis $\Gamma$
       contains an odd cycle and we are done), but $\displaystyle
       \bigcap_{i=1}^mF_i=\emptyset$. 

       So for each $j \in \{1,\ldots,m\}$, we can find a vertex $x_j$
       such that $x_j \in F_i \iff j \neq i$. Therefore we have a
       sequence of vertices $x_1,\ldots,x_m$ such that for each $i$:
       $$\{x_1, \ldots, \hat{x_i},\ldots,x_m\} \subseteq F_i \mbox{
       and } x_i \notin F_i.$$

      Now consider three facets $F_1,F_2,F_3$ of $\Gamma$. If
      $\{F_1,F_2,F_3\}$ is not a cycle, since it has length 3, it must
      be a tree; therefore it has a leaf, say $F_1$, and a joint, say
      $F_2$. It follows that $F_1 \cap F_3 \subseteq F_2$. But then it
      follows that $x_2 \in F_1 \cap F_3 \subseteq F_2$, which is a
      contradiction.  
     \end{proof}

We now concentrate on $L(\D)$ and its relation to $\D$.

\begin{lemma}\label{l:linegraph-subcollection} If  $\D$ is a facet complex, 
for every induced subgraph $G$ of $L(\D)$ there is a subset $\Gamma
\subseteq \D$ such that $G=L(\Gamma)$.
\end{lemma}

   \begin{proof} Let $G$ be an induced subgraph of $L(\D)$. Then if $V$ is 
   the vertex set of $\D$, and $W$ is the vertex set of $G$,
   $W\subseteq V$ and for each $x,y \in W$, $x$ and $y$ are connected
   by an edge in $G$ if and only if they are connected by an edge in
   $L(\D)$. This means that, if $F_1,\ldots,F_m$ are the facets of
   $\D$ corresponding to the vertices in $W$, and
   $\Gamma=\{F_1,\ldots,F_m\}$, then $G$ is precisely $L(\Gamma)$.
     \end{proof}

\begin{lemma}\label{l:cycle-graph} If $\D$ is a facet complex and $L(\D)$
is a cycle of length $\ell >3$, then $\D$ is a cycle of length $\ell$. 
\end{lemma}

    \begin{proof} Suppose $L(\D)$ is the cycle $$\{w_1,w_2\} 
     \sim_{L(\D)} \{w_2,w_3\} \sim_{L(\D)} \cdots \sim_{L(\D)}
     \{w_{\ell-1},w_{\ell}\}\sim_{L(\D)} \{w_{\ell},w_1\},$$ where
     each vertex $w_i$ of $L(\D)$ corresponds to a facet $F_i$ of
     $\D$. Since $w_i$ is only adjacent to $w_{i-1}$ and $w_{i+1}$
     (mod $\ell$), it follows that
     $$F_i\cap F_j \neq \emptyset \iff j=i-1,i, i+1 \mbox{ (mod
     $\ell$)}$$ which implies that $\D=\{F_1,\ldots,F_{\ell}\}$ where
     $$F_1\sim_{\D} F_2\sim_{\D} \cdots \sim _{\D} F_{\ell}\sim_{\D}
     F_1.$$ Moreover, since $\ell >3$, we have $\displaystyle
     \bigcap_{i=1}^{\ell}
     F_i=\emptyset$. Theorem~\ref{t:cycle-structure} now implies that
     $\D$ is a cycle of length $\ell$.
    \end{proof}

\begin{proposition}[The line graph of an odd-cycle-free complex]\label{p:odd-cycle-free-perfect} 
If $\D$ is an odd-cycle-free facet complex, then $L(\D)$ is either
perfect, or contains the complement of a 7-cycle as an induced
subgraph.
\end{proposition}

     \begin{proof}  Suppose $L(\D)$ is not perfect, and let $G$ be a 
      minimal imperfect induced subgraph of $L(\D)$. By
      Lemma~\ref{l:linegraph-subcollection}, for some subset $\Gamma$
      of $\D$, $G=L(\Gamma)$. By Theorem~\ref{t:SPGT}, $G$ is either
      an odd cycle of length $\geq 5$, or the complement of one.  If
      $G$ is an odd cycle, then so is $\Gamma$ by
      Lemma~\ref{l:cycle-graph}, and therefore $\D$ is not
      odd-cycle-free and we are done.

      So assume that $G$ is the complement of an odd cycle of length
      $\ell \geq 5$. We consider two cases.

       \begin{enumerate} 
    
        \item $\ell=5$: Since the complement of a 5-cycle is a 5-cycle, 
         it immediately follows from the discussions above that $\Gamma$ is
         a cycle of length 5, and hence $\D$ is not odd-cycle-free.

        \item $\ell \geq 9$: We show that $\Gamma$ contains a cycle of
        length 3. 
         
         Let $G=\overline{C_{\ell}}$, where $C_{\ell}$ is the $\ell$-cycle 
                  $$\{w_1,w_2\} \sim_{C_{\ell}} \{w_2,w_3\}
          \sim_{C_{\ell}} \cdots \sim_{C_{\ell}}
          \{w_{\ell-1},w_{\ell}\}\sim_{C_{\ell}} \{w_{\ell},w_1\},$$
          and a vertex $w_i$ of $G$ corresponds to a facet $F_i$ of
          $\Gamma$. This means that $F_i\cap F_j \neq \emptyset$
          unless $j=i-1,i+1$ (mod $\ell$). With this indexing,
          consider the subset $\Gamma'= \{F_1,F_4,F_7\}$ of
          $\Gamma$. Clearly all three facets of $\Gamma'$ have
          nonempty pairwise intersections:
         $$F_1 \cap F_4 \neq \emptyset, \ F_1 \cap F_7 \neq
           \emptyset,\ F_4 \cap F_7 \neq \emptyset.$$ Suppose
           $\Gamma'$ is not a cycle. Since $\Gamma'$ has only three
           facets it must be a tree and should therefore have a leaf,
           say $F_1$, and a joint, say $F_4$. So
           \begin{eqnarray}\label{e:joint} \emptyset \neq F_1\cap F_7
           \subseteq F_4.\end{eqnarray}

          Now consider the subset $\Gamma''=\{F_1,F_3,F_7\}$ of
          $\D$. We know that
          \begin{eqnarray}\label{e:intersections}F_1 \cap F_3 \neq
          \emptyset, \ F_1 \cap F_7 \neq \emptyset,\ F_3 \cap F_7 \neq
          \emptyset.\end{eqnarray} If $F_1\cap F_3 \cap F_7 \neq
          \emptyset$, then from (\ref{e:joint}) we see that $F_3 \cap
          F_4 \neq \emptyset$, which is a contradiction.  Therefore
          $F_1\cap F_3 \cap F_7 =\emptyset$, which along with the
          properties in (\ref{e:intersections}) and
          Lemma~\ref{l:3-cycle} implies that $\Gamma''$ is not a tree,
          so it must be a cycle.
   
          We can make similar arguments if $F_1$ or $F_7$ are joints
          of $\Gamma'$: if $F_1$ is a joint, then we can show that
          $\Gamma''=\{F_2,F_4,F_7\}$ is a cycle, and if $F_7$ is a
          joint, then $\Gamma''=\{F_1,F_4,F_6\}$ is a cycle. So we have
          shown that either $\Gamma'$ is a 3-cycle, or one can form
          another 3-cycle $\Gamma''$ in $\D$. Either way, $\D$
          contains an odd cycle, and is therefore not odd-cycle-free.
       \end{enumerate}
      \end{proof}

\idiot{We can argue in the one-before-last paragraph that since
odd-cycle-free complexes satisfy Helly property, $F_1\cap F_3 \cap F_7
\neq \emptyset$, and so the ``if'' part of that argument is not
necessary. On the other hand, we don't really need a strong property
like Helly here, since we get by easily without. So I left it as it is.}

Propositions~\ref{p:odd-cycle-free-perfect}
and~\ref{p:odd-cycle-free-Helly}, along with
Theorem~\ref{t:Helly-perfect-normal} immediately imply
Theorem~\ref{t:odd-cycle-free-normal}. Putting it all together, we
have shown that

\begin{theorem}[Odd-Cycle-Free complexes that satisfy K\"onig]\label{t:mainthm}
 If $\D$ is a facet complex that is odd-cycle-free and $L(\D)$ does
not contain the complement of a 7-cycle as an induced subgraph, then
every subset of $\D$ satisfies the K\"onig property.
\end{theorem}

      \begin{proof} The statement follows immediately from
        Theorem~\ref{t:odd-cycle-free-normal} along with
        Theorem~\ref{t:normal-konig}.
      \end{proof}

\subsection*{Are theses conditions necessary for satisfying K\"onig?}

A natural question is whether the conditions in
Theorem~\ref{t:mainthm} are necessary for a facet complex whose every
subset satisfies the K\"onig property. The answer in general is
negative. In this section, we explore various properties and examples
related to this issue.

The first observation is that not even all odd cycles fail the K\"onig
property. Indeed, if the cycle $\D$ (or in fact any complex) is a
cone, in the sense that all facets share a vertex, then it always
satisfies the K\"onig property with $\al(\D)=\be(\D)=1$.  

But if we eliminate the case of cones, all remaining odd cycles fail
the K\"onig property.

\begin{lemma}[Odd cycles that fail K\"onig]\label{l:odd-cycle-fail-konig}
 Suppose the facet complex $\D=\{F_1,\ldots,F_{2k+1}\}$ is a cycle of
odd length such that $\displaystyle
\bigcap_{i=1}^{2k+1}F_i=\emptyset$.  Then $\D$ fails the K\"onig
property.
\end{lemma}

    \begin{proof} Suppose without loss of generality that $\D$ can be 
     written as
     $$F_1\sim_\D F_2 \sim_\D \cdots \sim_\D F_{2k+1}\sim_\D F_1.$$
     Then a maximal independent set of facets of $\D$ can have at most
     $k$ facets; say $B=\{F_1,F_3,\ldots, F_{2k-1}\}$ is such a set,
     and by symmetry, all maximal independent sets will consist of
     alternating facets, and will have cardinality $k$. Hence
     $\be(\D)=k$.

     But we need at least $k+1$ vertices to cover $\D$. To see this,
     suppose that $\D$ has a vertex cover
     $A=\{x_1,\ldots,x_k\}$. Since $B$ is an independent set, we can
     without loss of generality assume that
     $$x_1 \in F_1, x_2 \in F_3, \ldots, x_i \in F_{2i-1}, \ldots, x_k
     \in F_{2k-1}.$$ 

     The other facets $F_2,F_4,\ldots,F_{2k},F_{2k+1}$ have to also be
     covered by the vertices in $A$. Since $F_{2k+1} \cap G=\emptyset$
     for all $G\in B$ except for $G=F_1$, we must have $x_1 \in
     F_{2k+1}$. Working our way forward in the cycle, and using the
     same argument, we get $$x_2 \in F_2, x_3 \in F_4, \ldots, x_i \in
     F_{2i-2},\ldots, x_k \in F_{2k-2}.$$ But we have still not
     covered the facet $F_{2k}$, who is forced to share a vertex of
     $A$ from one of its two neighbors: either $x_1 \in F_{2k}$ or
     $x_{k}\in F_{2k}$. Neither is possible as $F_{2k} \cap F_1=F_{2k}
     \cap F_{2k-2}=\emptyset$, and so $A$ cannot be a vertex cover.

    Adding a vertex of $F_{2k}$ solves this problem though, so
    $\al(\D)=k+1$, and hence $\D$ fails the K\"onig property.
    \end{proof}

The previous lemma then brings us to the question: can we replace the
condition ``odd-cycle-free'' with ``odd-hole-free'' (where an
\emph{odd hole} is referring to an odd cycle that is not a cone) in
the statement of Theorem~\ref{t:mainthm}? The answer is again
negative, as clarified by the example below.

\begin{example} The hollow tetrahedron $\D$ pictured below is 
odd-hole-free (but it does contain four 3-cycles which are
cones). However it fails the K\"onig property, since $\be(\D)=1$, but
$\al(\D)=2$.  Similar examples in higher dimensions can be constructed.

\begin{center}
 \epsfig{file=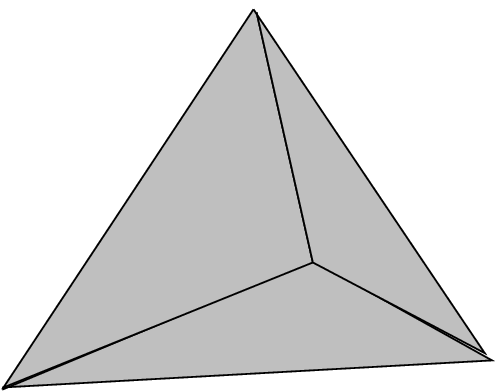, height=1in}
\end{center}
\end{example}

We next focus on the second condition in the statement of
Theorem~\ref{t:mainthm}, which turns out to be inductively necessary
for satisfying the K\"onig property.

\begin{lemma}\label{l:7-cycle-no-konig} Let $\D$ be a facet complex
such that $L(\D)$ is the complement of a 7-cycle. Then $\D$ fails the
K\"onig property.\end{lemma}

      \begin{proof} Suppose $L(\D)=\overline{C_7}$ where $C_7$ is a 
           7-cycle, and let $\D=\{F_1,\ldots,F_7\}$ such that the
          vertices of $C_7$ correspond to the facets $F_1,\ldots,F_7$
          in that order; in other words, $F_1$ intersects all other
          facets but $F_2$ and $F_7$, and so on (see
          Figure~\ref{f:compl-7-cycle}).

        \begin{figure}
         \begin{center}
          \epsfig{file=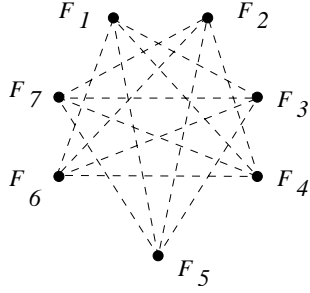, height=1.5in}
          \caption{Complement of a 7-cycle}\label{f:compl-7-cycle}
         \end{center}
	\end{figure}

       Let $B$ be a maximal independent set of facets, and assume $F_1
       \in B$. Then, since $F_1$ intersects all facets but $F_2$ and
       $F_7$, $B$ can contain one of $F_2$ and $F_7$ (but not both,
       since they intersect). So $|B|=2$. The same argument holds if
       $B$ contains any other facet than $F_1$, so we conclude that
       $\be(\D)=2$.

      Now suppose $\D$ has a vertex cover of cardinality 2, say
      $A=\{x,y\}$. Then each facet of $\D$ must contain one of $x$ and
      $y$. Without loss of generality, suppose $x \in F_1$. Since each
      facet does not intersect the next one in the sequence
      $F_1,F_2,\ldots,F_7, F_1$, we have
      $$x \in F_1 \Longrightarrow
        y \in F_2 \Longrightarrow
        x \in F_3 \Longrightarrow
        y \in F_4 \Longrightarrow
        x \in F_5 \Longrightarrow
        y \in F_6 \Longrightarrow
        x \in F_7.$$

      But now $x\in F_1 \cap F_7 =\emptyset$, which is a contradiction. So 
      $\al(\D)\geq 3$, and hence $\D$ does not satisfy the K\"onig property.
      \end{proof}

\begin{corollary} If every subset of a facet complex $\D$ satisfies 
the K\"onig property, then $L(\D)$ cannot contain the complement of a
7-cycle as an induced subgraph.
\end{corollary}

\begin{remark}[The case of the complement of a 7-cycle] 
As suggested above, if $L(\D)$ contains the complement of a 7-cycle as an
induced subgraph, $\D$ may fail the K\"onig property, even though it
may be odd-cycle-free. For example, consider the complex $\D$ on seven
vertices $x_1,\ldots,x_7$: $\D=\{F_1,\ldots,F_7\}$ where
$F_1=\{x_1,x_2,x_3\}$, $F_6=\{x_2,x_3,x_4\}$, $F_4=\{x_3,x_4,x_5\}$,
$F_2=\{x_4,x_5,x_6\}$, $F_7=\{x_5,x_6,x_7\}$, $F_5=\{x_6,x_7,x_1\}$,
$F_3=\{x_7,x_1,x_2\}$.

The graph $L(\D)$ is the complement of a 7-cycle (the labels of the
facets correspond to those in Figure~\ref{f:compl-7-cycle}). One can
verify that $\D$ contains no 3, 5, or 7-cycles, so it is
odd-cycle-free. However by Lemma~\ref{l:7-cycle-no-konig}, the facet
complex $\D$ fails the K\"onig property; indeed $\al(\D)=3$ but
$\be(\D)=2$.

On the other hand, it is easy to expand $\D$ to get another complex
$\Gamma$, such that $L(\Gamma)$ does contain the complement of a
7-cycle as an induced subgraph, and $\Gamma$ satisfies the K\"onig
property. For example, consider $\Gamma=\{G,F_1',F_2,\ldots,F_7\}$,
where $F_2,\ldots, F_7$ are the same facets as above, and we introduce
two new vertices $u,v$ to build the new facets
$F_1'=\{u,x_1,x_2,x_3\}$, and $G=\{u,v\}$.

The set $B=\{G,F_2,F_3\}$ is a maximal independent set of facets, so
$\be(\Gamma)=3$. Also, we can find a vertex covering
$A=\{u,x_4,x_7\}$, which implies that $\al(\D)=3$.

Note, however, that $\Gamma$ does not satisfy the K\"onig property
``inductively'': it contains a subset $\{F_1',F_2,\ldots,F_7\}$
that fails the K\"onig property by Lemma~\ref{l:7-cycle-no-konig}.
\end{remark}

\section{Balanced complexes are odd-cycle-free}\label{s:hyper-cycles}

The notion of a cycle has already been defined in hypergraph theory,
and is much more general than our definition of a cycle
(see~\cite{B2}, or Chapter 5 of~\cite{B1}). To keep the terminologies
separate, in this paper we refer to the traditional hypergraph cycles
as \emph{hyper-cycles}. In particular, hypergraphs that do not contain
hyper-cycles of odd length are known to satisfy the K\"onig
property. In this section, we introduce this class of hypergraphs and
show that hypergraphs not containing odd hyper-cycles are
odd-cycle-free, and their line graphs cannot contain the complement of
a 7-cycle as an induced subgraph.

\begin{definition}[Hyper-cycle~\cite{B1,B2}]\label{d:hyper-cycle}
 Let $\HH$ be a hypergraph on vertex set $V$. A \emph{hyper-cycle} of
length $\ell$ ($\ell\geq 2$), is a sequence
$(x_1,F_1,x_2,F_2,\ldots,x_{\ell},F_{\ell},x_1)$ where the $x_i$ are distinct
vertices and the $F_i$ are distinct facets of $\HH$, and moreover
$x_i,x_{i+1} \in F_i$ (mod $\ell$) for all $i$.
\end{definition}

\begin{definition}[Balanced hypergraph (\cite{B1,B2})] A hypergraph is said to be 
\emph{balanced} if every odd hyper-cycle has an edge containing three
vertices of the cycle.
\end{definition}

 Herzog, Hibi, Trung and Zheng~\cite{HHTZ} called a hyper-cycle a
 \emph{special cycle} if, with notation as in
 Definition~\ref{d:hyper-cycle}, for all $i$ we have $x_i \in F_j$ if
 and only if $j=i-1,i$ (mod $\ell$). In other words, if each vertex
 $x_i$ of the hyper-cycle appears in exactly two facets, the
 hyper-cycle is a special cycle. So a balanced hypergraph is one that
 does not contain any special cycle of odd length. Special cycles have
 also been called \emph{strong} cycles in the literature.

It is easy to see that a cycle $\D$ defined as
$$F_1\sim_{\D} F_2\sim_{\D} \cdots \sim _{\D} F_{\ell}\sim_{\D} F_1$$
produces a hyper-cycle; just pick any vertex $x_i \in F_i \cap
F_{i+1}$ (mod $\ell$), $\D$ produces a hyper-cycle, or in fact a
special cycle, of the same length $\ell$
$$(x_1,F_1,x_2,F_2,\ldots,x_{\ell},F_{\ell},x_1).$$

It follows that a balanced simple hypergraph is odd-cycle-free.  The
converse, however, is not true.

\begin{example}[Not all odd-cycle-free complexes are balanced]\label{e:ocf-not-balanced}
 Consider the complex $\D$ below, which is odd-cycle-free, as the only
 cycle is the 4-cycle $\{H_1,H_2,H_4,H_5\}$. But $\D$ is not balanced,
 as all of $\D$ forms the special 5-cycle
 $$(x_6,H_1,x_2,H_2,x_3,H_3,x_4,H_4,x_4,H_5,x_6).$$

\begin{center}
 \epsfig{file=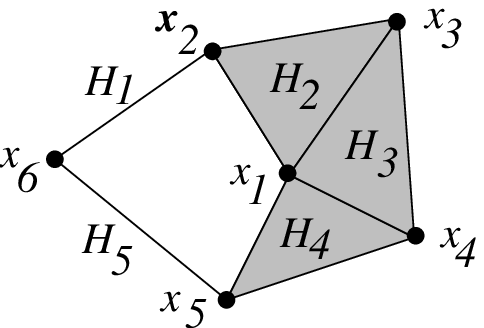, height=1.2in}
\end{center}
\end{example}

The complex in Example~\ref{e:ocf-not-balanced} is an example of how
our main result (Theorem~\ref{t:mainthm}) generalizes the fact that
balanced complexes satisfy the K\"onig property. In fact, we can show
the following.

\begin{proposition}\label{p:balanced-ocf} Let $\D$ be a balanced complex. 
Then $\D$ is odd-cycle-free and $L(\D)$ does not contain the
complement of a 7-cycle as an induced subgraph. \end{proposition}

    \begin{proof} Let $\D$ be balanced. We have already shown that
     $\D$ is odd-cycle-free. Suppose that $L(\D)$ contains the
     complement of a 7-cycle as an induced subgraph. 

          By lemmas~\ref{l:linegraph-subcollection}
          and~\ref{l:cycle-graph}, $\D$ contains a subset $\Gamma$
          whose line graph is the complement of a 7-cycle, we write
          $L(\Gamma)=\overline{C_7}$. Suppose
          $\Gamma=\{F_1,\ldots,F_7\}$ such that the vertices of the
          7-cycle $C_7$ correspond to the facets $F_1,\ldots,F_7$ in
          that order; in other words, $F_1$ intersects all other
          facets but $F_2$ and $F_7$, and so on (see
          Figure~\ref{f:compl-7-cycle}).

          We claim that we can find vertices $x_1 \in F_1\cap F_5$, $x_2 \in 
          F_5\cap F_7$, $x_3 \in F_2\cap F_7$, $x_4 \in F_2\cap F_6$, and $x_5
          \in F_1\cap F_6$, such that
          \begin{eqnarray}\label{e:5-cycle}
          (x_5,F_1,x_1,F_5,x_2,F_7,x_3,F_2,x_4,F_6,x_5)
          \end{eqnarray}
          is a special cycle of length 5, and if not, $\D$ contains a
          special cycle of length 3.

          The main obstacle in making the hyper-cycle in
          (\ref{e:5-cycle}) a special 5-cycle, is finding appropriate
          choices for the vertices $x_1,\ldots,x_5$. Suppose these
          choices are not possible, in which case at least one of the
          following statements hold:

           \begin{enumerate}
   
           \item $F_1 \cap F_5 \subseteq F_2 \cup F_6 \cup F_7$. 

                 This is not possible, since we know that $F_1 \cap
                 F_2=F_1\cap F_7 =\emptyset$, and $F_5 \cap
                 F_6=\emptyset$. Since $F_1 \cap F_5 \neq \emptyset$,
                 one can choose $x_1 \in F_1\cap F_5$ such that $x_1
                 \notin F_2 \cup F_6 \cup F_7$.

           \item $F_5 \cap F_7 \subseteq F_1 \cup F_2 \cup F_6
           \Longrightarrow F_5 \cap F_7 \subseteq F_2$ (Since $F_7
           \cap F_1=F_7 \cap F_6=\emptyset$).

             In this case, consider the facet complex
             $\{F_3,F_5,F_7\}$. Then, since $F_2 \cap F_3=\emptyset$
             and $F_5 \cap F_7 \subseteq F_2$, we have $F_3 \cap F_5
             \cap F_7=\emptyset$. Lemma~\ref{l:3-cycle} now implies
             that $\{F_3,F_5,F_7\}$ is a 3-cycle, and hence can be
             written as a special 3-cycle.

           \item $F_2 \cap F_7 \subseteq F_1 \cup F_5 \cup F_6
           \Longrightarrow F_2 \cap F_7 \subseteq F_5$ (Since $F_7
           \cap F_1=F_7 \cap F_6=\emptyset$).

            Similar to Case 2. it follows that $\{F_2,F_4,F_7\}$ is a
            (special) 3-cycle.

           \item $F_2 \cap F_6 \subseteq F_1 \cup F_5 \cup F_7$.

                 Fails with argument similar to Case 1. So one can
                 choose $x_4 \in F_2\cap F_6$ such that $x_4 \notin
                 F_1 \cup F_5 \cup F_7$.

           \item $F_1 \cap F_6 \subseteq F_2 \cup F_5 \cup F_7$.

                 Fails with argument similar to Case 1. So one can
                 choose $x_5 \in F_1\cap F_6$ such that $x_5 \notin
                 F_2 \cup F_5 \cup F_7$.

	   \end{enumerate}

           So we have shown that either there are vertices
           $x_1,\ldots,x_5$ such that the sequence in
           (\ref{e:5-cycle}) is a special 5-cycle, or otherwise,
           either cases 2. or 3. above would hold, in which case $\D$
           would contain a (special) 3-cycle. Either way, $\D$ is not
           balanced.
       \end{proof}

As a result, we have another proof to the following known fact
(see~\cite{B1,B2}).

\begin{corollary}[Balanced complexes satisfy K\"onig] If $\D$ is a balanced 
facet complex, then all subsets of $\D$ satisfy the K\"onig property.
\end{corollary}

In fact, a stronger version of the above statement was proved for
balanced hypergraphs by Berge and Las Vergnas; see page 178
of~\cite{B1}.

In closing, we would like remark that odd-cycle-free facet complexes,
unlike balanced ones, are not necessarily Mengerian. As indicated
below, this fact has implications for the algebraic properties of the
facet ideals of such complexes.

\begin{remark}[Odd-Cycle-Free complexes are not Mengerian] 
The facet complex in Example~\ref{e:ocf-not-balanced} is an example of
an odd-cycle-free complex which is not Mengerian. To see this, let
$M$ be the incidence matrix of $\D$, where the rows correspond to the facets and the columns to the vertices of $\D$:
$$M=\left (
\begin{tabular}{llllll}
0&1&0&0&0&1\\
1&1&1&0&0&0\\
1&0&1&1&0&0\\
1&0&0&1&1&0\\
0&0&0&0&1&1\\
\end{tabular}
\right )
$$ and pick the vector ${\bf c}=(2, 1, 1, 1, 2, 2)$. 

For $\D$ to be Mengerian, we need 
$$\mbox{min}\{{\bf a}.{\bf c} \ |\ {\bf a} \in ({\mathbb{N}}\cup \{0\})^6,\ M {\bf a } \geq {\bf 1}\}
=
\mbox{max}\{{\bf b}. {\bf 1} \ |\ {\bf b} \in ({\mathbb{N}}\cup \{0\})^5,\ M^T{\bf b } \leq {\bf c}\}.
$$

One can check that with this value of ${\bf c}$, the left-hand minimum
value is $4$, but the right-hand maximum value is 3.

This implies that odd-cycle-free complexes, unlike bipartite graphs or
balanced complexes, do not necessarily have normally-torsion-free
facet ideals (see \cite{HHTZ, SVV}).

Computational evidence using the computer algebra softwares
\emph{Normaliz}~\cite{BK} and \emph{Singular}~\cite{GPS} indicates,
however, that these ideals may still be normal. It would be of great
interest to know whether odd-cycle-free complexes provide a new class
of normal ideals; this would generalize results in~\cite{F1, HHTZ,
SVV}.
\end{remark}


\end{document}